\documentstyle[11pt]{article}
\input{psfig.sty}
\input{amssymb.sty}
\normalbaselines
\begin{document}
\bibliographystyle{unsrt}
\newtheorem{thm}{{\sc Theorem}}[section]
\newtheorem{lma}{{\sc Lemma}}[section]
\newtheorem{prp}{{\sc Proposition}}[section]
\newtheorem{cor}{{\sc Corollary}}
\newtheorem{lemma}{{Lemma}}
\newtheorem{theorem}{{Theorem}}

\newfont{\msbm}{msbm10 scaled\magstep1}
\newfont{\eusm}{eusm10 scaled\magstep1}

\def\cpx{{{\Bbb C}}}
\def\spec{{\mbox{spec}}}
\def\Hom{{\mbox{Hom}}}
\def\itg{{{\Bbb Z}}}
\def\real{{{\Bbb R}}}
\def\rat{{{\Bbb Q}}}
\def\eM{{{I(\widehat{o})}^{\bot}}}
\def\Io{{I(\widehat{o})}}
\def\bas{{\mbox{\eusm B}}}
\def\wt{{\mbox{wt}}}
\def\QT{{\mbox{\goth{Q}}}}
\def\PT{{\mbox{\goth{P}}}}
\def\TT{{\mbox{\goth{T}}}}

\def\bea*{\begin{eqnarray*}}
\def\eea*{\end{eqnarray*}}
\def\ba{\begin{array}}
\def\ea{\end{array}}
\count1=1
\def\be{\ifnum \count1=0 $$ \else \begin{equation}\fi}
\def\ee{\ifnum\count1=0 $$ \else \end{equation}\fi}
\def\ele(#1){\ifnum\count1=0 \eqno({\bf #1}) $$ \else \label{#1}\end{equation}\fi}
\def\req(#1){\ifnum\count1=0 {\bf #1}\else \ref{#1}\fi}
\def\bea(#1){\ifnum \count1=0   $$ \begin{array}{#1}
\else \begin{equation} \begin{array}{#1} \fi}
\def\eea{\ifnum \count1=0 \end{array} $$
\else  \end{array}\end{equation}\fi}
\def\elea(#1){\ifnum \count1=0 \end{array}\label{#1}\eqno({\bf #1}) $$
\else\end{array}\label{#1}\end{equation}\fi}
\def\cit(#1){
\ifnum\count1=0 {\bf #1} \cite{#1} \else 
\cite{#1}\fi}
\def\bibit(#1){\ifnum\count1=0 \bibitem{#1} [#1    ] \else \bibitem{#1}\fi}
\def\ds{\displaystyle}
\def\hb{\hfill\break}
\def\comment#1{\hb {***** {\em #1} *****}\hb }
\def\Pr{{{\Bbb P}}}

\newcommand{\TZ}{\hbox{\bf T}}
\newcommand{\MZ}{\hbox{\bf M}}
\newcommand{\NZ}{\hbox{\bf N}}
\def\ZZ{{\itg}}
\def\RZ{{\real}}
\def\CZ{{\cpx}}
\def\PZ{{\Pr}}
\def\QZ{{\rat}}
\newcommand{\HZ}{\hbox{\bf H}}
\newcommand{\EZ}{\hbox{\bf E}}
\newcommand{\GZ}{\,\hbox{\bf G}}
\font\germ=eufm10
\def\goth#1{\hbox{\germ #1}}
\newtheorem{pf}{Proof}
\renewcommand{\thepf}{}
\vbox{\vspace{38mm}}
\begin{center}
{\LARGE \bf Orbifolds and Finite Group  
Representations }\\[5mm]
Li Chiang 
\\{\it Institute of Mathematics \\ Academia Sinica \\ 
Taipei , Taiwan \\
(e-mail:chiangl@gate.sinica.edu.tw)}
\\[5mm]
Shi-shyr Roan
\\{\it Institute of Mathematics \\ Academia
Sinica \\  Taipei , Taiwan \\ (e-mail:
maroan@ccvax.sinica.edu.tw)} \\[5mm]
\end{center}

\begin{abstract} 
We present our recent understanding on 
resolutions of Gorenstein orbifolds, which
involves the finite group representation theory. We shall concern only the quotient singularity of hypersurface type. The abelian group $A_r(n)$ for  
 $A$-type hypersurface quotient singularity of
dimension $n$ is introduced. For $n=4$,
the structure of Hilbert scheme of group
orbits  and crepant resolutions of $A_r(4)$-singularity are obtained. The
flop procedure of 4-folds is explicitly
constructed through the process.  
\par \vspace{5mm} \noindent
1991 MSC 14J, 14M, 20C. 
\par \vspace{2mm} \noindent

\end{abstract}

\vfill
\eject

\section{Introduction}
It is well known that the theory of
``minimal" resolutions of singularity of algebraic
(or analytical) varieties differs significantly
when the  (complex) dimension of the variety is
larger than two. As the prime achievement in
algebraic geometry of the 1980-s, the minimal
model program in the 3-dimensional birational
geometry carried out by S. Mori and others has
provided an effective tool
for the study of algebraic 3-folds,
(see
\cite{Mo} and references therein). Meanwhile,
Gorenstein quotient singularities in dimension 3
has attracted considerable 
interests among geometers
due to the development of string theory, by which
the orbifold Euler characteristic of an orbifold
was proposed as the vacuum description of
models built upon the quotient of a
manifold 
\cite{DHVW}. The consistency of physical
theory then demanded the existence of crepant
resolutions which are compatible with the orbifold
Euler characteristic. The complete
mathematical justification of the conjecture was obtained
in the mid-90s (see \cite{Rtop} and references
therein). However, due to the computational
nature of methods in the proof, the
qualitative understanding of the these crepant
resolutions has still been lacking on certain
aspects from a mathematical viewpoint. Until
very recently, by the 
development of Hilbert scheme of a finite group
$G$-orbits, initiated by Nakamura and others
with the result obtained 
\cite{BKR, GNS, INj, IN, Na}, it strongly indicates a
promising role of the finite group in problems of
resolutions of quotient singularities.
In particular a plausible method has been suggested on the study of
geometry of orbifolds through the group
representation theory. It has been known that
McKay correspondence \cite{Mc} between
representations of Kleinian groups and affine
$A$-$D$-$E$ root diagrams has revealed a profound
geometrical structure on the minimal
resolution of the quotient surface singularity
(see e.g.,
\cite{GV}). A similar connection between the
finite group and general quotient singularity
theories would be expected. Yet, the interest
of this interplay of geometry and group
representations would not only aim on the research
of crepant resolutions, but also on its own
right, due to possible implications on
understanding some certain special kind of group
representations by engaging the rich algebraic
geometry techniques. 

In this article, we shall study problems
related to the  crepant resolutions of quotient
singularities of higher dimension
$n$, (mainly for $n \geq 4$).  Due to the
many complicated exceptional cases of the
problem,  we shall restrict ourselves here only
on those of the hypersurface singularity type.
The purpose of this paper is to present certain 
primitive results of our first attempt on the study of the higher
dimensional hypersurface orbifolds  under the
principle of ``geometrization" of finite group
representations.  We shall give a brief account of
the progress recently made. The main issue we
deal with in this work will be the higher
dimensional generalization of the $A$-type
Kleinian surface singularity, the
$A_r(n)$-hypersurface singularity of dimension
$n$ (see (\ref{An2}) of the content). For $n=4$,
we are able to determine the detailed structure of
$A_r(4)$-Hilbert scheme, and its relation
with crepant resolutions of $\CZ^4/A_r(4)$.
In the process,
an explicit ``flop" construction of $4$-folds among
different crepant resolutions is found.
In this article, 
we shall  only sketch 
the main ideas behind the proof of these results,
referring the reader to our forecoming paper 
\cite{CR} for a more  complete  description of
methods and  arguments  used.

This paper is organized as follows. 
In Sect. 2, we shall give a brief introduction
of the general scheme of engaging finite group
representations in the birational geometry of
orbifolds. Its connection with the Hilbert scheme
of $G$-orbits for a finite linear group $G$ on
$\CZ^n$, ${\rm Hilb}^G(\CZ^n)$, introduced in
\cite{IN}, will be explained in Sect. 3.
In Sect. 4, we first review certain basic
facts in toric geometry, which will be presented 
in the form most suitable for our goal, then focusing
the case on $A_r(n)$-singularity. For $n=3$, we will give a
thorough discussion on the explicit toric
structure of
${\rm Hilb}^{A_r(3)}(\CZ^3)$ as an illustration
of  the general result obtained by
Nakamura on abelian group $G$ in \cite{Na}\footnote{For a finite group in ${\rm SL}_3(\CZ)$, a recent result in \cite{BKR} has shown that ${\rm Hilb}^G(\CZ^3)$ is a crepant resolution of $\CZ^3/G$.}. In Sect. 5, we deal with  a
special case of $4$-dimensional orbifold with $G=
A_1(4)$, and derive the detailed structure of
${\rm Hilb}^G(\CZ^4)$. Its relation with the
crepant resolutions of $\CZ^4/G$ is given, so
is the ``flop" relation among crepant
resolutions. In Sect. 6, we describe the result
of 
$G= A_r(4)$ for the arbitrary $r$ , then end with
some  concluding remarks.

{\bf Notations. } To present 
our work, we prepare some notations. In this
paper, by an orbifold we shall always mean the
orbit space for  a finite
group action on a smooth complex 
manifold. For a finite group
$G$, we denote 
\begin{eqnarray*}
{\rm Irr}(G) &= & \{ \rho : G \longrightarrow
{\rm GL}( V_{\rho}) \ \  {\rm an
\ irreducible
\ representative \ of \ } G \} .
\end{eqnarray*}
The trivial representation of $G$ will be denoted
by ${\bf 1}$. For a $G$-module $W$,
i.e., a $G$-linear representation on a vector
space 
$W$, one has the canonical irreducible
decomposition:
$$
W = \bigoplus_{\rho \in {\rm Irr}(G)} W_{\rho} 
$$
where $W_{\rho}$ is a $G$-submodule of
$W$, $G$-isomorphic to $V_{\rho} \otimes
W_{\rho}^0$ for some  trivial $G$-module 
$W_{\rho}^0$. The vector space $W_{\rho}$ will
be called the
$\rho$-factor of the $G$-module $W$.

For an analytic variety $X$, we shall not distinct the notions of 
vector bundle and locally free ${\cal
O}_X$-sheaf over $X$. For  a vector bundle $V$
over $X$, an automorphisms of $V$ means a  linear
automorphism with the identity on 
$X$. If the bundle
$V$ is acted by a group $G$ as bundle
automorphisms, we  shall call
$V$ a
$G$-bundle.

\section{Representation Theory in Algebraic
Geometry of Orbifolds } In this paper, $G$ will
always denote a finite (non-trivial) subgroup of
${\rm GL}_n(\CZ)$ for
$n \geq 2$, and $S_G : = \CZ^n/G$ with the
canonical projection,
$$
\pi_G: \CZ^n \ \longrightarrow \ S_G \ ,
$$
and $o:= \pi_G(0) \in S_G$. 
When $G$ is a subgroup of ${\rm SL}_n(\CZ)$, which
will be our main concern later in
this paper,  
$G$ acts on $\CZ^n$ freely outside a finite
collection of linear  subspaces with codimension
$\geq 2$. Then the orbifold $S_G$ has a non-empty
singular set, ${\rm Sing}(S_G)$, of codimension
$\geq 2$, in fact,
$o \in {\rm Sing}(S_G)$. 

For $G$ in ${\rm GL}_n(\CZ)$, 
$S_G$ is a singular variety in general.
By a birational morphism of a variety
over $S_G$, we shall always  mean a proper
birational morphism $\sigma$ from variety $X$ to
$S_G$  
which defines a  biregular map between $X
\setminus \sigma^{-1}({\rm Sing}(S_G))$ and $ S_G
\setminus {\rm Sing}(S_G) $,
\be
\sigma : X \longrightarrow  S_G \ .
\ele(sigma)
One has the commutative diagram,
\bea(cll)
X \times_{S_G} \CZ^n & \longrightarrow & \CZ^n \\ 
\ \ \downarrow \pi & & \downarrow \pi_G  \\
X  &\stackrel{\sigma}{\longrightarrow} & S_G  \ .
\elea(Xdia)
Denote ${\cal F}_X$ the coherent ${\cal O}_X$-sheaf over $X$ obtained by 
the push-forward of the structure sheaf of $X \times_{S_G} \CZ^n$,
$$
{\cal F}_X : = \pi_* {\cal O}_{X \times_{S_G}\CZ^n} \ .
$$
The sheaf ${\cal F}_X$ 
has the following 
functorial property, namely for $X,
X'$ birational  over $S_G$ with the commutative
diagram,
$$
\begin{array}{lll}
\ \ X'  &\stackrel{\sigma'}{\longrightarrow} & S_G \\ 
\mu \downarrow & &  || \\
\ \ X &\stackrel{\sigma}{\longrightarrow} & S_G  
\end{array}
$$
one has an canonical morphism, 
$ \mu^*{\cal F}_X \longrightarrow {\cal F}_{X'}$. 
In particular, with the morphism (\req(sigma)) 
we have the ${\cal O}_X$-morphism,
$$
\sigma^* ( {\pi_G}_*{\cal O}_{\CZ^n} )
\longrightarrow {\cal F}_X  \ .
$$
Furthermore, all the morphisms in the above are
compatible with the natural $G$-structure on
${\cal F}_X$ induced from the
$G$-action on
$\CZ^n$ via (\req(Xdia)). One has the
canonical $G$-decomposition of ${\cal
F}_X$,
$$
{\cal F}_X = \bigoplus_{\rho \in {\rm Irr}(G)}
({\cal F}_X)_\rho \ ,
$$
where $({\cal F}_X)_\rho$ is the $\rho$-factor of 
${\cal F}_X$, and it is a coherent ${\cal O}_X$-sheaf over $X$. The geometrical fiber of ${\cal
F}_X, ({\cal F}_X)_\rho$ over an element $x$
of $X$ are defined by
$$
{\cal F}_{X, x} = k(x) \bigotimes_{{\cal O}_X} {\cal F}_X \ , \ 
({\cal F}_X)_{\rho, x} = k(x) \bigotimes_{{\cal
O}_X} ({\cal F}_X)_\rho \ , 
$$
where $k (x) ( := {\cal O}_{X, x}/{\cal M}_{x})$
is the residue field  at $x$. Over 
$X - \sigma^{-1}({\rm Sing}(S_G))$, 
${\cal F}_X$ 
is a vector bundle of rank $|G|$  with the
regular $G$-representation on each geometric
fiber. Hence $({\cal F}_X)_\rho$ is a vector
bundle over 
$X - \sigma^{-1}({\rm Sing}(S_G))$ of the 
rank equal to
${\rm dim}. V_\rho$. For $x \in X$, there
exists a
$G$-invariant ideal $I(x)$ in $\CZ[Z] (:=
\CZ[Z_1, \cdots, Z_n ])$ such that 
the following relation holds,
\begin{eqnarray}
{\cal F}_{X, x} = k (x)\bigotimes_{{\cal O}_{S_G} }
{\cal O}_{\CZ^n} 
(x) \simeq \CZ[Z]/I(x) \ .
\label{I(x)} 
\end{eqnarray}
We have $({\cal F}_X)_{\rho, x}
\simeq (\CZ[Z]/I(x))_\rho$. 
The vector spaces
$\CZ[Z]/I(x)$ form a family of finite
dimensional $G$-modules parameterized by
$x \in X$, which are equivalent
to the regular representation for elements
outside 
$\sigma^{-1}({\rm Sing}(S_G))$.

Set  $X=S_G$ in (\ref{I(x)}). For 
$s \in S_G$, there is a
$G$-invariant ideal $I(s)$ of $\CZ[Z]$; in
fact, $I(s)$ is the ideal in $\CZ[Z]$
generated by the $G$-invariant polynomials
vanishing at $\sigma^{-1}(s)$.
Let $\widetilde{I}(s)$  be the ideal of $\CZ[Z]$
consisting of all polynomials in $\CZ[Z]$
vanishing at
$\sigma^{-1}(s)$. Then $\widetilde{I}(s)$ is an
$G$-invariant ideal, and we have
$$
\widetilde{I}(s) \supset I(s)  \ .
$$
In
particular, for 
$s= o$, we have 
$$  
\widetilde{I}(o) = \CZ[Z]_0 \ , \ \ 
I(o) = \CZ[Z]^G_0 \CZ[Z] 
$$
where the subscript $0$ means the maximal ideal
with polynomials vanishing at the
origin.  For a birational variety $X$ over $S_G$
via $\sigma$ in (\req(sigma)), the
following relations of $G$-invariant ideals of
$\CZ[Z]$  hold:
$$
\widetilde{I}(s) \supset I(x) \supset
I(s) \ , \ \ x \in X \ , \ s = \sigma(x) \ . 
$$
A certain connection exists
between algebraic geometry and $G$-modules through the
variety $X$. 
For $x \in X$, 
there is a direct sum $G$-decomposition of
$\CZ[Z]$,
$$
\CZ[Z] = I(x)^{\bot} \oplus I(x) \ .
$$
Here $I(x)^{\bot}$ is a finite dimensional
subspace of $\CZ[Z]$ which is $G$-isomorphic to
$\CZ[Z]/I(x)$. Similarly we have the
$G$-decomposition of $\CZ[Z]$ for $s=\sigma(x)
\in S_G$,
$$
\CZ[Z] = I(s)^{\bot} \oplus
I(s) \ ,
\ \ \CZ[Z] = \widetilde{I}(s)^{\bot} \oplus
\widetilde{I}(s)
$$ 
such that the following relations hold for the
finite-dimensional $G$-modules, 
$$
\widetilde{I}(s)^{\bot} \subset I(x)^{\bot}
\subset I(s)^{\bot}
\ . 
$$
Consider the canonical
$G$-decomposition of $I(x)^{\bot}$,
$$
I(x)^{\bot} = \bigoplus_{\rho \in {\rm Irr}(G)}
I(x)^{\bot}_\rho \ .
$$
Note that $I(x)^{\bot}_\rho $ is isomorphic to
a positive finite copies of $V_{\rho}$. Then the 
affine structure of $X$ near $x$ is determined by 
the $\CZ$-algebra generated by all  the
$G$-invariant rational function $f(Z)$ such that
$f(Z) I(x)^{\bot}_\rho
\subset I(x)$ for some
$\rho$.

\section{Hilbert Scheme of  Finite Group
Orbits}

Among the varieties       
$X$ birational over $S_G$ with ${\cal F}_X$ a
vector bundle,  there exists a 
minimal object, called the $G$-Hilbert scheme in 
\cite{IN, Na},   
$$
\sigma_{\rm Hilb} : {\rm Hilb}^G(\CZ^n)
\longrightarrow  S_G \  .
$$ 
For another $X$, the map (\req(sigma)) can be factored through 
a birational morphism $\lambda$ from $X$ onto
${\rm Hilb}^G(\CZ^n)$ via $\sigma_{\rm Hilb}$,
$$
\lambda : X \longrightarrow {\rm Hilb}^G(\CZ^n) \
.
$$ 
In fact, the ideal $I (x)$,$x \in X$, of
(\ref{I(x)}) are with  the co-length $|G|$,  which
gives rise to the above map 
$\lambda$ of $X$ to ${\rm Hilb}^G(\CZ^n)$. We
shall denote $X_G$ as the normal variety
over $S_G$ defined by 
\begin{eqnarray*}
X_G := {\rm normalization \ of \ } {\rm
Hilb}^G(\CZ^n) \ , \ \ \sigma_G: X_G
\longrightarrow S_G \ .
\end{eqnarray*}
By the fact that  every biregular
automorphism of
$S_G$ can always be lifted to one on 
${\rm Hilb}^G(\CZ^n)$, hence on $X_G$, one has the
following result.
\begin{lemma}\label{lem:Aut(SG)}
Denote ${\rm Aut}(S_G)$ the group of biregular automorphisms 
of $S_G$. Then ${\rm Hilb}^G(\CZ^n)$ and $ X_G$ are 
varieties over $S_G$ with the ${\rm
Aut}(S_G)$-equivariant covering morphisms.
\end{lemma}
By the definition of ${\rm Hilb}^G(\CZ^n)$, an element $p$ of ${\rm Hilb}^G(\CZ^n)$ represents a 
$G$-invariant ideal $I(p)$ of $\CZ[Z]$ of
co-length $|G|$. The fiber of the vector
bundle ${\cal F}_{{\rm Hilb}^G(\CZ^n)}$ over $p$
can be identified with he regular $G$-representation space 
$\CZ[Z]/I(p)$.
Our study will mainly concern on the
relation of crepant resolutions of $S_G$ and
${\rm Hilb}^G(\CZ^n)$. For this purpose we will
assume for the rest of the
paper that the group $G$ is
a subgroup of
${\rm SL}_n(\CZ)$:
$$
G \subset {\rm SL}_n(\CZ) \ ,
$$ 
which is the same to say that $S_G$ has the Gorenstein quotient
singularity.
For $n=2$, these groups were classified
by F. Klein into $A$-$D$-$E$ types \cite{Kl}, the
singularities are called Kleinian singularities. The
minimal resolution $\widehat{S}_G$ of $S_G$ has
the trivial canonical bundle ( i.e., crepant), by
\cite{HH}. In \cite{IN, Na},  Y. Ito and I.
Nakamura showed that ${\rm Hilb}^G(\CZ^2)$ is
equal to the minimal resolution 
$\widehat{S}_G$. For $n=3$, it has been known
that there exist crepant resolutions
for a 3-dimensional Gorenstein orbifold (see
\cite{Rtop} and references therein). Two different crepant
resolutions of the same orbifolds are connected by
a sequences of flop processes (see e.g.,
\cite{R90}). It was expected that   ${\rm
Hilb}^G(\CZ^3)$ is one of those
crepant resolutions. The assertion has been
confirmed in the abelian group case in
\cite{Na}, and in general by \cite{BKR}. 

For the motivation of our later study on the
higher dimensional singularities, we now 
illustrate the relation between 
$G$-Hilbert scheme and the minimal resolution in
dimension 2, i.e., surface singularities.
For the rest of this section, we are going to
describe the structure of
${\rm Hilb}^G(\CZ^2)$ for the
$A$-type Kleinian group,
\begin{eqnarray}
G = A_r := \{  \left( \begin{array}{cc}
\epsilon &   0  \\
0 & \epsilon^{-1}
\end{array} \right) \ | \ \epsilon^{r+1}=1 \} \ ,
\ \ r \geq 1 \ . \label{An2}
\end{eqnarray}
The affine ring of $\CZ^2$ is $\CZ[Z] (=\CZ[Z_1,
Z_2])$ and $G$-invariant polynomials is the
algebra  generated by 
$$
Y_1:= Z_1^{r+1} \ , \ Y_2 := Z_2^{r+1} \ , \
X := Z_1Z_2  \ .
$$
Hence the ideal $I(o)$ in $\CZ[Z]$ for $o \in
S_G$ is equal to
$\langle Z_1^{r+1}, Z_1^{r+1}, Z_1Z_2 \rangle$,
and
$ Z_1^k, Z_2^k , (  0 \leq k \leq r ),$ form a
basis of the $G$-module $I(o)^{\bot}$. For a
non-trivial character $\rho$,
$I(o)^{\bot}_{\rho}$ is of dimension 2, in fact,
it has a basis consisting of a pair, $Z_1^k,
Z_2^{r+1-k}$, for some $k$. With the method of
continued fraction \cite{HH}, it is known that
the minimal resolution
$\widehat{S}_G$ of
$S_G$ has the trivial canonical bundle with an open affine cover $\{ {\cal U}_k
\}_{k=0}^r $, where the
coordinates $(u_k, v_k)$ of ${\cal U}_k$ is expressed by 
$$
{\cal U}_k \simeq \CZ^2 \ni (u_k, v_k) =
(Z_1^{k+1}Z_2^{-r+k} , Z_1^{-k}Z_2^{r+1-k} ) \ .
$$
Denote $\hat{o}_k$ the element in $\widehat{S}_G$
with the coordinate $u_k=v_k=0$. The exceptional
divisor in $\widehat{S}_G$
is $E_1+ \cdots +E_r$ where $E_j$ is  
a rational $(-2)$-curve joining
$\hat{o}_{j-1}$ and $\hat{o}_j$. The configuration can be realized in
the following tree diagram:
$$
\put(-105, 0){\line(3,-1){35}}
\put(-80, -12){\line(3,1){35}}
\put(-60, 0){\line(3,-1){35}}
\put(-40, -12){\line(3,1){35}}
\put(-20, 0){\line(3,-1){35}}
\put(2, -12){\line(3,1){35}}
\put(20, 0){\line(3,-1){35}}
\put(40, -12){\line(3,1){40}}
\put( 67, 2){\line(3,-1){40}}
\put(100, -20){\shortstack{$\hat{o}_0$ }}
\put(70, -10){\shortstack{$\hat{o}_1$ }}
\put(45, -20){\shortstack{$\hat{o}_2$ }}
\put(-80, -20){\shortstack{$\hat{o}_{r-1}$ }}
\put(-108, -12){\shortstack{$\hat{o}_r$ }}
\put(-15, -13){\shortstack{$\hat{o}_j$ }}
\put(90, -2){\shortstack{$E_1$ }}
\put(50, -2){\shortstack{$E_2$ }}
\put(0, -3){\shortstack{$E_j$ }}
\put(-90, -2){\shortstack{$E_r$ }}
\put(-180, -40){\shortstack{[Fig. 1] \
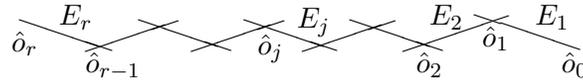
Exceptional curve configuration in the
minimal resolution of $\CZ^2/A_r$. }}
$$ 
It is easy to see that the ideal $I({\hat{o}_k})$ is
given by
$$
I({\hat{o}_k}) = \langle Z_1^{k+1}, Z_2^{r+1-k},
Z_1Z_2 \rangle  , 
$$
hence the $G$-module $\CZ[Z]/I({\hat{o}_k})$ is
the regular representation isomorphic the
following one,
\be
I({\hat{o}_k})^\bot = \CZ + 
\sum_{i=1}^k \CZ Z_1^i + \sum_{j =1}^{r-k} \CZ
Z_2^j \ .
\ele(Iok)
One can represent monomials in the above
expression as the ones with $\bullet$ in the
following picture:
$$
\put(-80, -80){\line(1, 0){107}}
\put(-80, 25){\line(0, -1){105}}
\put(19, -83) { \shortstack{$\times$}}
\put(4, -83) { \shortstack{$\times$}}
\put(-11, -83) { \shortstack{$\times$}}
\put(-26, -83) { \shortstack{$\times$}}
\put(-41, -82) { \shortstack{$\bullet$}}
\put(-86, -82) { \shortstack{$\bullet$}}
\put(-71, -82) { \shortstack{$\bullet$}}
\put(-56, -82) { \shortstack{$\bullet$}}
\put(-86, -67) { \shortstack{$\bullet$}}
\put(-86, -52) { \shortstack{$\bullet$}}
\put(-86, -37) { \shortstack{$\bullet$}}
\put(-86, -24) { \shortstack{$\bullet$}}
\put(-88, -9) { \shortstack{$\times$}}
\put(-88, 6) { \shortstack{$\times$}}
\put(-88, 21) { \shortstack{$\times$}}
\put(-93, 21){ \shortstack{$r$}}
\put(-115, -24){ \shortstack{$r-k$}}
\put(-93, -69){ \shortstack{$1$}}
\put(20, -90){ \shortstack{$r$}}
\put(-40, -90){ \shortstack{$k$}}
\put(-70, -90){ \shortstack{$1$}}
\put(-210, -105){ \shortstack{[Fig. 2] \ 
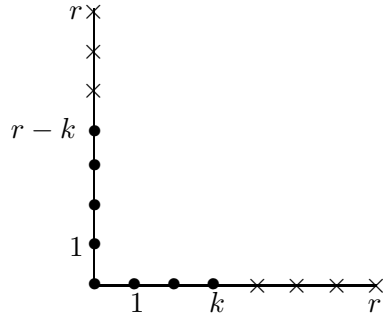
Representatives of $I({\hat{o}_k})^\bot
(\simeq \CZ[Z]/I({\hat{o}_k}))$ for the minimal
resolution of $\CZ^2/A_r$.}}
$$
For $x \in {\cal U}_k$,  the ideal $I(x)$ has the
expression
$$
I(x) = \langle Z_1^{k+1}- \alpha Z_2^{r-k} ,
Z_2^{r+1-k} - \beta Z_1^k , Z_1Z_2 - \alpha \beta
\rangle \ , \ \ \alpha, \beta \in \CZ .
$$
The classes in 
$\CZ[Z]/I(x)$ represented by monomials in
(\req(Iok)) still form a basis, hence give rise
to a local frame of the vector bundle ${\cal
F}_{\widehat{S}_G}$ over ${\cal U}_k$. The divisor   
$E_{k+1}$ is defined by $\beta=0$, and its element approaches to 
$\hat{o}_{k+1}$ as $\alpha$ tends to
infinity.

\section{Abelian Orbifolds and Toric Geometry}
In this section we discuss the abelian group case of $G$ in the previous section using methods in toric geometry. We shall consider $G$ as a subgroup of the diagonal group of
${\rm GL}_n(\CZ)$, denoted by $T_0 := \CZ^{* n}$,
and regard 
$\CZ^n$ as the partial compactification of $T_0$,
\begin{eqnarray*}
G \subset T_0 
\subset \CZ^n \ .
\end{eqnarray*}
Define the $n$-torus $T$ with the toric embedding
in $S_G  (= \CZ^n/G)$ by 
$$
T := T_0/G \ , \ \ \ T \subset S_G \ .
$$
Techniques in toric geometry rely 
on lattices of one-parameter subgroups, 
characters of $T_0, T$, 
\begin{eqnarray*}
N (:= {\rm Hom}(\CZ^*, T )) &\supset& N_0 (:= {\rm
Hom}(\CZ^*, T_0 )) \ , \\
M (:= {\rm Hom}( T, \CZ^* )) &\subset& M_0 (:=
{\rm Hom}( T_0, \CZ^* )) \ .
\end{eqnarray*}
For our convenience, we shall make the following
identification of $N_0, N$ with lattices in
$\RZ^n$. An element $x$ in $\RZ^n$ has the 
coordinates $x_i$ with respective to the 
standard basis $(e^1, \cdots , e^n )$:
\[ 
x = \sum_{i=1}^n x_ie^i  \in \RZ^n . 
\]
Then 
$$
N_0  = \ZZ^n (:= {\rm exp}^{-1}(1)) \ , \ \ 
N = {\rm exp}^{-1}( G ) \ , 
$$
where ${\rm exp}: {\RZ}^n \longrightarrow T_0 $
is defined by ${\rm exp}(x) = \sum_{i} e^{2\pi
\sqrt{-1}x_i}e^i$. Note that $G \simeq N/N_0$.
The dual lattice $M_0$ of $N_0$ is the standard one in the dual space $\RZ^{n *}$, and we shall identify it with the group of monomials of $Z_1,
\ldots, Z_n$ via the correspondence:
$$
I = \sum_{s=1}^n i^se_s \in M_0 \ \ 
\longleftrightarrow \ \ Z^I = \prod_{s=1}^n
Z_s^{i_s}
$$  
The dual lattice $M$ of $N$ is the sublattice of $M_0$
corresponding to the set of $G$-invariant monomials.

Over the $T$-space $S_G$, we now consider only those 
varieties 
$X$ which is normal
and birational over $S_G$ with a $T$-structure, hence as it has been known, are presented by certain 
combinatorial data by the toric method \cite{Dan, KKMS, O}. Note that by Lemma (\ref{lem:Aut(SG)}), $X_G$ is a toric variety over $S_G$. 
In general, a toric variety over $S_G$ is
described by  a fan
$\Sigma = \{ \sigma_{\alpha} \ | \ \sigma \in I
\}$ whose support equals to  the first
quadrant of
$\RZ^n$, i.e., 
a rational convex cone decomposition of the first
quadrant of
$\RZ^n$. Equivalently, it is determined by the
intersection of the fan and the simplex
$\triangle$ where 
$$
\bigtriangleup : = \{ x \in {\RZ}^n | \sum_{i}
x_i = 1, x_j \geq 0 \ \ \forall \ j \} .
$$
The data in $\bigtriangleup $ is given by
$\Lambda  = \{
\bigtriangleup_{\alpha} \ | \ \alpha \in I \} $, 
where $
\bigtriangleup_{\alpha}:= \sigma_{\alpha} \bigcap 
\bigtriangleup$. The  $\bigtriangleup_{\alpha} $s
form a decomposition of
$\bigtriangleup$ by convex subsets, having the
vertices in $\bigtriangleup \cap
\QZ^n$. Note that for $\sigma_{\alpha}=
\{ \vec{0} \}$, we have
$\bigtriangleup_{\alpha}= \emptyset$. We shall
call  
$\Lambda $ a
rational polytope decomposition of
$\bigtriangleup$, and denote the corresponding
toric variety by $X_\Lambda$. We call $\Lambda$
an integral polytope decomposition of
$\bigtriangleup$ if all the vertices of 
$\Lambda$ are in $N$. For a 
rational polytope decomposition $\Lambda$ of
$\bigtriangleup$, we define $\Lambda(i):= \{
\triangle_{\alpha} \in \Lambda \ | \ {\rm
dim}(\triangle_{\alpha}) = i \}$ for $ -1 \leq i
\leq n-1$, (here ${\rm dim}( \emptyset ):= -1$).
Then
$T$-orbits in
$X_\Lambda$ are parameterized by
$\bigsqcup_{i=-1}^{n-1} \Lambda(i)$.  In fact,
for each
$\bigtriangleup_{\alpha} \in \Lambda(i)$, there
associates a 
$T$-orbit of the
dimension $n-1-i$, denoted by
${\rm orb}(\bigtriangleup_{\alpha})$. 
A toric divisor in $X_\Lambda$ is the closure
of a $n-1$ dimensional orbit, denoted by $D_v =
\overline{{\rm orb}(v)}$ for $v \in 
\Lambda (0)$. The canonical sheaf of 
$X_{\Lambda}$ has the expression in terms of toric divisors (see, e.g. \cite{KKMS}),
\begin{equation}
\omega_{X_\Lambda} = {\cal O}_{X_\Lambda}(
\sum_{v \in \Lambda(0) } (m_v-1)D_v ) \ , \ 
\label{K}
\end{equation}
where $m_v$ is the  positive integer 
such that
$m_v v$ is a primitive element of 
$N$. In particular, the crepant property of $X_\Lambda$, i.e.  
$\omega_{X_\Lambda}= {\cal O}_{X_\Lambda}$, is
given by the integral condition of $\Lambda$.
The non-singular criterion of  $X_{\Lambda}$
is the simplicial  decomposition of $\Lambda$ together 
with the multiplicity one property, i.e., for
each $\Lambda_{\alpha} \in \triangle(n-1)$,
the elements $m_v v$, $v \in \Lambda_{\alpha} \cap \Lambda(0)$, 
form a $\ZZ$-basis of $N$. 
 The
following results are known for toric variety
over $S_G$  (see e.g. \cite{R89} and
references therein):

(1) The Euler number of $X_\Lambda$ is given by
$\chi ( X_\Lambda ) = | \Lambda (n-1) |$.

(2) For a rational polytope
decomposition $\Lambda$ of $\triangle$, any
two of the following three
conditions implies the third
one:
\begin{eqnarray*} 
X_\Lambda : \mbox{ non-singular
} , \ \ \ \omega_{X_\Lambda} = {\cal O}_{X_\Lambda} \
, \ \ \ \chi ( X_\Lambda ) = | G | . 
\end{eqnarray*}
It is easy to see that the following result holds for the sheaf ${\cal F}_{X_\Lambda}$.
\begin{lemma} 
Let $\Lambda$ be a
rational polytope decomposition of
$\bigtriangleup$, and $x_0$ be the zero-dimensional
toric orbit in 
$X_\Lambda$ corresponding to an element   
$\bigtriangleup_{\alpha_0}$ in
$\Lambda (n-1) $. Let $Z^{I^{(j)}},
1
\leq j \leq N$, be a
finite collection of monomials whose classes
generate the $G$-module $\CZ[Z]/I(x_0)$. Then the
classes of $Z^{I^{(j)}}$s also generate
$\CZ[Z]/I(y)$ for $y \in {\rm
orb}(\bigtriangleup_{\beta})$ with $
\bigtriangleup_{\beta} \subseteq 
\bigtriangleup_{\alpha_0}$.
\end{lemma}

Define 
\begin{eqnarray}
A_r(n) := \{ g \in {\rm SL}_n({\CZ}) \ | \ g:
\mbox{ diagonal} \ ,\  g^{r+1}=1  \} \ , \ \ r \geq 1 \ .
\end{eqnarray}
Note that  the above group for $n=2$ is the same as $A_r$ in (\ref{An2}). 
For a general $n$, 
$A_r(n)$-invariant polynomials in $\CZ[Z]$
are generated by the following
$(n+1)$ ones:
$$
X := \prod_{i=1}^n Z_i \ , \ \ \   Y_j :=
Z_j^{r+1}
\ (j=1,
\ldots, n) \ .
$$
This implies that $S_{A_r(n)}$ is  the 
singular hypersurface in $\CZ^{n+1}$,
$$
S_{A_r(n)} = \{ (x, y_1, \cdots, y_n ) \in
\CZ^{n+1} \ | \  x^{r+1} = y_1
\cdots y_n \  \} \ .
$$
For the rest of this
paper, we will conduct the discussion of abelian
orbifolds mainly on the group $A_r(n)$.  The ideal $I(o)$ of 
$\CZ[Z]$ associated to the element $o \in
S_{A_r(n)}$ is given by 
$$
I(o) = <Z_1^{r+1}, \ldots, Z_n^{r+1}, Z_1 \cdots
Z_n > \ , 
$$
hence
$$
I(o)^{\bot} = \bigoplus \{ \CZ Z^I \ | \ I=(i^1,
\ldots, i^n) \ , \ \ 0 \leq i^j \leq r \ ,
\prod_{j=1}^n i^j=0 
\} \ .
$$
For ${
\bf 1} \neq \rho \in {\rm
Irr}(A_n(r))$, the dimension of $I(o)^{\bot}_\rho$ is always greater than one.
In fact, one can describe  explicitly a set of  
monomial generators of 
$I(o)^{\bot}_\rho$. For example, say
$I(o)^{\bot}_\rho$ containing an element $Z^I$ with 
$I=(i^1, \ldots, i^n), i^1=0$ and $ i^s \leq
i^{s+1}$, then $I(o)^{\bot}_\rho$ is generated by
$Z^K$s with 
$K=(k^1,  \ldots, k^n)$ given by   
\be
k^s= \left\{ \begin{array}{ll} 
r+1-i^j + i^s & {\rm if} \  i^s <i^j \ , \\
i^s - i^j & {\rm otherwise} \ ,
\end{array} \right. 
\ele(Ioorg)
here  $j$ runs through $1$ to $ n$. 
In particular for $r=1$, the dimension of 
$I(o)^{\bot}_\rho$ is equal to $2$ for $\rho 
\neq {\bf 1}$, with a basis consisting of $Z^I, Z^{I'}$
whose indices satisfy the relations, $0 \leq i^s,
i^{s'}
\leq 1 , i^s+i^{s'}=1$ for $1 \leq s \leq n$. 

For $n=3$, by the general result of Nakamura on an abelian group $G$
(Theorem 4.2 in \cite{Na}), ${\rm
Hilb}^{A_r(3)}(\CZ^3)$ is a crepant toric
variety. To illustrate this fact, we give here a
direct derivation of the result by working
on the explicitly described toric variety.
\par \vspace{.1in} \noindent
{\bf Example}. It is easy to see that $\triangle
\cap N$ consists of the following elements, 
$$
v^{m_1,m_2,m_3} := \frac{1}{r+1}\sum_{j=1}^3
m_je^j
\ , \ \ \ \ \  m_j \in \ZZ_{\geq 0} \ , \ \
\sum_{j=1}^3 m_j = r+1 \ .
$$
Denote $\Xi$ the simplicial decomposition of
$\triangle$ obtained by drawing the three lines 
parallel to the edges of $\triangle$ through each element of $\triangle \cap N$. The 2-simplexes in $\Xi(2)$ consists of two types
of triangles, $\triangle^{m_1, m_2, m_3}_u, 
\triangle^{m_1, m_2, m_3}_d$:
\begin{eqnarray*}
\triangle^{m_1, m_2, m_3}_u &= \langle
v^{m_1, m_2, m_3}, \ v^{m_1-1, m_2 +1,m_3}, \
v^{m_1-1, m_2, m_3+1} \rangle  \ ,  \\
\triangle^{m_1, m_2, m_3}_d &= 
\langle
v^{m_1, m_2, m_3}, \ v^{m_1-1, m_2+1,m_3}, \
v^{m_1, m_2+1, m_3-1} \rangle \ ,
\end{eqnarray*}
the vertices of each one form a $\ZZ$-basis of $N$.

\begin{figure}[hh]
\begin{center}
\mbox{\psfig{figure=triadecm.ps}}
\mbox{\psfig{figure=localdec.ps}}
\vspace{0.5cm}

[Fig. 3] The left one is the simplicial decomposition $\Xi$ of $\triangle$ for $r=4$. \\ 
The right one is the union of $\triangle^{m,l,k}_u$ and $\triangle^{m,l,k}_d$.
   
\end{center}
\end{figure}
\noindent
Therefore the $X_{\Xi}$ is a smooth toric variety. The $0$-dimensional
(toric) orbits corresponding to $\triangle^{m_1, m_2, m_3}_u, 
\triangle^{m_1, m_2, m_3}_d$ are denoted by
$x^{m_1, m_2, m_3}_u$, $ x^{m_1, m_2, m_3}_d$
respectively.  The affine coordinate system
centered at $x^{m_1, m_2, m_3}_u$ is given by
$$
U_1 = \frac{Z_1^{r+2-m_1}}{(Z_2Z_3)^{m_1-1}} , \ \ \ 
U_2=\frac{Z_2^{r+1-m_2}}{(Z_1Z_3)^{m_2}} , \ \ \ U_3=
\frac{Z_3^{r+1-m_3}}{(Z_1Z_2)^{m_3}}
\ .
$$
For $y$ in $X_{\Xi}$ with the value 
$U_j$ equal to $\alpha_j$, by computation
the ideal
$I(y)$ has the following generators,
$$
\begin{array}{lll}
Z_1^{ r+2-m_1}-\alpha_1(Z_2Z_3)^{m_1-1}, & 
Z_2^{ r+1-m_2}-\alpha_2(Z_3Z_1)^{m_2}, & 
Z_3^{ r+1-m_3}-\alpha_3(Z_1Z_2)^{m_3}, \\
(Z_1Z_2)^{m_3+1} -\alpha_1\alpha_2 Z_3^{ r-m_3}, &
(Z_2Z_3)^{m_1} -\alpha_2\alpha_3 Z_1^{ r+1-m_1},&
(Z_3Z_1)^{m_2+1} -\alpha_3\alpha_1 Z_2^{ r-m_2},\\
 Z_1Z_2Z_3-\alpha_1\alpha_2\alpha_3 \ . & \ &
\end{array}
$$
For the element $ x^{m_1, m_2, m_3}_u$, i.e., $\alpha_i=0$
for all $i$, it is easy to see that there are
$(r+1)^2$ monomials not in 
$I(x^{m_1, m_2, m_3}_u)$, which form a basis of $I(x^{m_1, m_2,
m_3}_u)^\bot$, and  give rise to a
basis of $\CZ[Z]/I(y)$ for $y$ in the affine 
neighborhood of $ x^{m_1, m_2, m_3}_u$. 
The same argument applies to the affine
neighborhood near $x^{m_1, m_2, m_3}_u$ with  the coordinate functions,
$$
V_1= \frac{(Z_2Z_3)^{m_1}}{Z_1^{r+1-m_1}} , \ \ \
V_2= \frac{(Z_1Z_3)^{m_2+1}}{Z_2^{r-m_2}}, \ \ \ 
V_3= \frac{(Z_1Z_2)^{m_3}}{Z_3^{r+1-m_3}} \ ,
$$ 
hence the description of ideals for $y$
in
$X_{\Xi}$ with
$V_j=\beta_j$,  
$$
\begin{array}{lll}
Z_1^{r+2-m_1}-\beta_2\beta_3(Z_2Z_3)^{m_1-1}, &
Z_2^{r+1-m_2}-\beta_3\beta_1(Z_3Z_1)^{m_2},  &
Z_3^{r+2-m_3}-\beta_1\beta_2(Z_1Z_2)^{m_3-1}, \\ 
(Z_1Z_2)^{m_3}-\beta_3Z_3^{r+1-m_3} , &
(Z_2Z_3)^{m_1}-\beta_1Z_1^{r+1-m_1} ,&
(Z_3Z_1)^{m_2+1}-\beta_2Z_2^{r-m_2} ,
\\
Z_1Z_2Z_3 -\beta_1\beta_2\beta_3 \ . & &
\end{array}
$$
Therefore we have shown that $X_{\Xi}$ is birational 
over ${\rm Hilb}^G(\CZ^3)$. Now we
are going to show that they are in fact
the same. Let $x$ be an element in ${\rm Hilb}^G(\CZ^3)$
represented by  an monomial ideal $J=I(x)$, (i.e,
with  a set of generators composed of monomials). Then the regular
$G$-module $J^{\bot}$ is generated by $|G|$
monomials, and $x$ lies over the element $o$ of $S_G$,
equivalently,
$J$ contains the ideal $\CZ[Z]^G_0$. Denote $l_i$
the smallest non-negative integer such that
$Z_i^{l_i} \in J$,  $l_{i,j}$ the smallest non-negative integer with
$(Z_iZ_j)^{l_{i,j}}
\in J $ for $i \neq j$. Hence $1 \leq l_i \leq
r+1$, and $Z_i^{l_i-1} \in J^{\bot}$, which implies 
$(\frac{Z_1Z_2Z_3}{Z_i})^{r+2-l_i} \in J$. In
particular, $Z_1^{l_1-1} \in J^\bot$ and
$(Z_2Z_3)^{r+2-l_1} \in J$. 
By the description (\req(Ioorg)) for $I(o)^\bot$,
$Z_1^{l_1}$ is the only
monomial in the basis of $I(o)^\bot$ for the
corresponding character of $G$, the same for $ (Z_2Z_3)^{r+1-l_1}$. Hence
$(Z_2Z_3)^{r+2-l_1} \in J^\bot$, which implies
$l_1+l_{2,3}= r+2$. Similarly, we have 
$l_2+l_{1,3}=l_3+l_{1,2}= r+2$. Again by
(\req(Ioorg)), $Z_1^{l_1-1}Z_2^{l_{1,2}-1}$,
$Z_2^{r+1-l_1+l_{1,2}}Z_3^{r+2-l_1},
Z_1^{l_1-l_{1,2}}Z_3^{l_3}$ are the generators
of an eigenspace in $I(o)^\bot$. The latter two
are elements in $J$ by $(Z_2Z_3)^{r+2-l_1},
Z_3^{l_3} \in J$. Therefore
$Z_1^{l_1-1}Z_2^{l_{1,2}-1} \in J^\bot$.
Similarly, one has
$$
Z_2^{l_2-1}Z_1^{l_{1,2}-1}, 
Z_1^{l_1-1}Z_3^{l_{1,3}-1},
Z_3^{l_3-1}Z_1^{l_{1,3}-1},  
Z_2^{l_2-1}Z_3^{l_{2,3}-1}, 
Z_3^{l_3-1}Z_2^{l_{2,3}-1} \in J^\bot \ .
$$
This implies that 
$J$ is the ideal in $\CZ[Z]$ generated by 
$Z_i^{l_i} , (Z_jZ_k)^{l_{j,k}} , Z_1Z_2Z_3$, 
where $\{i, j, k\} =\{1,2,3\}$ with $ j<k $, and $l_i, l_{j,k}$ are integers satisfying the relations:
$l_i+l_{j,k}=r+2$ , $1 \leq l_i \leq r+1$.
Therefore the number of monomials in $J^\bot$ is
given by
$$
\begin{array}{l}
1 + \sum_{j=1}^3(l_j-1) +
\sum_{j<k}(l_j-1)(l_k-1)-
(l_j-l_{j,k})(l_k-l_{j,k}) \\ =
-(\sum_{j}l_j)^2+(4r+7)\sum_{j}l_j
-4(r+1)^2-6(r+1)-2 \ ,
\end{array}
$$
which is equal to $(r+1)^2$ by $J \in {\rm
Hilb}^G(\CZ^3)$. The only solutions for
$\sum_jl_j$ are $2(r+1)+1$, $2(r+1)+2 $. Hence the
ideal $J$ is characterized by integers $l_i$ between 1 and $r+1$  with the relation:
$$
\sum_{j=1}^3 l_j = 2(r+1)+1 , \ \ 2(r+1)+2 \ .
$$
If $\sum_{j=1}^3 l_j = 2(r+1)+1$, the ideal $J$
corresponds to the ideal of $x^{m_1,
m_2, m_3}_u$ in $X_{\Xi}$ with 
$$
l_1= r+2-m_1, \ l_2= r+1-m_2 , \ l_3= r+1-m_3 \ .
$$
For $\sum_{j=1}^3 l_j = 2(r+1)+2$, the ideal $J$
corresponds to the ideal for the element $x^{m_1,
m_2, m_3}_d$ in $X_{\Xi}$ with 
$$
l_1= r+2-m_1, \ l_2= r+1-m_2 , \ l_3= r+2-m_3 \ .
$$
Now by techniques of Gr\"{o}bner basis for ideals in
$\CZ[Z]$ (see, e.g., \cite{CLO}), 
given a monomial order, one have a monomial ideal $\mbox{lt}(J)$ (generated by the leading monomial of the Gr\"{o}bner basis of $J$) such that all monic monomials outside $\mbox{lt}(J)$ form a base of $\CZ[Z]/J$. Thus,
one sees that for an element in ${\rm Hilb}^G(\CZ^3)$, i.e., a $G$-invariant ideal $J'$ with the regular $G$-module $\CZ[Z]/J'$, there is a basis of 
$\CZ[Z]/J'$ represented by monomials in $J^\bot$ for some $J$ previously described. 
This shows ${\rm
Hilb}^G(\CZ^3) = X_{\Xi}$.
\hfill $\Box$ \par \vspace{.2in} \noindent
{\bf Remark.} As in the discussion of the case
for $n=2$ in Sect. 3, one can represent the
monomial basis elements of $I(x^{m_1,m_2,
m_3}_u)^\bot$ ,   $I(x^{m_1,m_2,
m_3}_d)^\bot$ in a pictorial way. For example,
the following ones are for $r=4$, and $m_1=2$, $m_2=1$, $m_3=2$:

\begin{figure}[hh]
\begin{center}
\mbox{\psfig{figure=plots1.ps}}
\mbox{\psfig{figure=plots2.ps}}
\vspace{0.5cm}

[Fig. 4] Graph representation of $I(x^{m_1,m_2,m_3}_u)^\bot$, $I(x^{m_1,m_2,m_3}_d)^\bot$ for $r=4$, $(m_1,m_2,m_3)=(2, 1, 2)$. An ``$\bullet$" indicates a monomial in $I^\bot$ while an ``$\times$" means one in $I$.  The difference between two graphs are marked by broken segments.
\end{center}
\end{figure}

\section{ ${\bf A_1(4)}$-Singularity and 
Flop  of 4-folds }
We now study the $A_r(n)$-singularity with $n
\geq 4$. For simplicity, let
us consider the case $r=1$, i.e.,  
$G= A_1(n)$, (indeed, no conceptual
difficulties arise for higher values of $r$).
The $N$-integral elements in $\triangle$ are as
follows: 
$$
\bigtriangleup \bigcap N = \{ e^j \ | \ 1 \leq j
\leq n \}  \bigcup \{ v^{i,j} \ | \ 1
\leq i < j \leq n \} \ , 
$$
where $ v^{i,j} :=
\frac{1}{2}(e^i + e^j )$ for $i \neq j$.
Other than the whole simplex $\triangle$, there
is only one integral polytope decomposition of 
$\triangle$ invariant under permutations of
coordinates, denoted by $\Xi$, which we now describe as follows. There are $n+1$ elements in
$\Xi (n-1)$:
$\bigtriangleup_{i} , 1 \leq i \leq n,$ together with 
$\Diamond$, where 
$\bigtriangleup_{i}
$ is the simplex generated by $e^i$ and
$v^{i,j}$ for $j \neq i$, and 
$\Diamond=$ the closure of $ \bigtriangleup
\setminus \bigcup_{i=1}^n
\bigtriangleup_{i}$. In fact, $\Diamond$ is the
convex hull spanned by all the $v^{i, j}$ for $i
\neq j$. The lower dimensional polytopes of
$\Xi $ are given by the faces of those in 
$\Xi (n-1)$. Then $X_{\Xi}$ has the trivial
canonical sheaf. However only for $n= 2, 3$,
$X_{\Xi}$  is a crepant resolution of
$S_{A_1(n)}$ (see, e.g.,
\cite{R89}). In fact, one has the following
result for higher $n$.
\begin{lemma}  \label{lem:SingXi}
For $n=4$, the toric variety $X_{\Xi}$ is smooth except one
isolated singularity, which is the 0-dimensional T-orbit
corresponding to $\Diamond$.
\end{lemma}
{\it Proof.} In general, for $n\ge 4$, it is easy to see that for each
$i$, the vertices of $\triangle_i$ form a
$\ZZ$-basis of $N$, for example, say $i=1$, it which follows
from $| A_1(n)| = 2^{n-1}$, and 
$$
{\rm det}(e^1, v^{1,2}, \cdots, v^{1,n} ) =
\frac{1}{2^{n-1}} \ .
$$
Hence $X_{\Xi}$ is nonsingular near the $T$-orbits associated to simplices in $\Delta_i$. As $\Diamond$ is 
not a simplex, ${\rm orb}(\Diamond)$ is always a singular point of $X_\Xi$. For $n=4$, the statement of smoothness of $X_\Xi$ except ${\rm orb}(\Diamond)$ follows from the fact that for $1\le i\le 4$, the vertices  $v^{i,j}$ ($j\ne i$) of $\Diamond$ together with $(1/2)\sum_{j=1}^4 e^j$, from an $N$-basis.
\hfill $\Box$ \par \vspace{.2in} \noindent
{\bf Remark.} (1) Denote $x_j := {\rm
orb}(\triangle_j) \in X_{\Xi}$ for $1 \leq j \leq
n$. The inverse of the matrix of vertices of
$\triangle_j$, 
$$
( v^{1,j} , \cdots, v^{j-1,j}, e^j,  v^{j+1,
j}, \cdots, v^{n, j})^{-1} \ ,
$$
gives rise to affine coordinates $(U_1,
\ldots, U_n)$ around $x_j$:  
$$
U_i = Z_i^2 \ \ (i \neq j) \ , \ \ \ U_j = 
\frac{Z_j}{Z_1
\cdots \check{Z}_j \cdots Z_n} \ .
$$
Hence
$$
I(x_j) = <Z_j , Z_i^2 , i \neq j  > + I_{A_1(n)} 
$$
with the regular $A_1(n)$-module isomorphism,
\bea(l)
\CZ[Z]/I(x_j) \simeq \bigoplus \{ \CZ Z^I \ | 
I=(i_1, \ldots, i_n) , \ i_j=0 , 
i_k =0, 
1 \ {\rm for } \ \ k
\neq j \} \ .
\elea(A1reg)

(2) We shall denote $x_{\Diamond}$ the element ${\rm
orb}(\Diamond)$ in $X_{\Xi}$, $x_{\Diamond}:= {\rm
orb}(\Diamond)$. The singular structure of
$x_{\Diamond}$ is determined by those $A_1(n)$-invariant polynomials, corresponding to 
$M$-integral elements in the cone dual to the
one generated by $\Diamond$ in $N_{\RZ}$. It is
easy to see that these polynomials are generated
by the following ones: 
$$
X_j := Z_j^2 \ , \ \ Y_j := \frac{Z_1 \cdots
\check{Z_j}
\cdots Z_n}{Z_j} \ .
$$
Hence we have
$$
I(x_{\Diamond}) = \langle Z_1 \cdots
\check{Z_j}
\cdots Z_n  \rangle_{1 \leq j \leq n}  + I_{A_1(n)} \ . 
$$
Note that for $n=3$, $Y_j$s form the
minimal generators for the invariant polynomials,
which implies the smoothness of $X_{\Xi}$. For $n \geq 4$,  
$x_{\Diamond}$ is an isolated singularity, but not of
the hypersurface type. For
$n=4$, the $X_j, Y_j (1 \leq j \leq 4)$ form a minimal set of 
generators of invariant polynomials, hence the
structure near $x_{\Diamond}$ in $X_{\Xi}$
is the 4-dimensional affine variety
in $\CZ^8$  defined by the relations:
\begin{eqnarray}
(x_i, y_i)_{1
\leq i \leq 4} \in \CZ^8 : \  x_iy_i= x_j y_j ,
\ \
\  x_ix_j = y_{i'}y_{j'} \ , \ 
\label{4sing}
\end{eqnarray}
where $i \neq j$ and $\{ i', j' \} $ is the
complimentary pair of $\{i, j\}$.
\hfill $\Box$ \par \vspace{.2in} \noindent

For the rest of this section, we shall consider
only the case $n=4$. We shall discuss the crepant resolutions of
$S_{A_1(4)}$, and its relation with 
${\rm Hilb}^{A_1(4)}(\CZ^4)$. Now the simplex 
$\bigtriangleup$ is a tetrahedron, and 
$\Diamond$ is an octahedron, on which the
symmetric group ${\goth S}_4$ acts as the
standard representation. The dual polygon of
$\Diamond$ is the cubic. Faces of
the octahedron $\Diamond$ are labeled by 
$F_j, F_j'$ for $1 \leq j \leq 4$, where 
$$
F_j = \Diamond \cap \bigtriangleup_j \ , \ \ \
F_j' = 
\{ \sum_{i=1}^4 x_je^j \in \Diamond \ | \ x_j =0
\} \ .
$$
The dual of $F_j$, $F_j'$ in the cubic are vertex, denoted by $\alpha_j$, $\alpha_j'$ as in [Fig. 5].

\begin{figure}[ht]
\begin{center}
\mbox{\psfig{figure=octeface.ps}}\hskip 0.5in
\mbox{\psfig{figure=cubic.ps}}

[Fig. 5] Dual pair of octahedron and cube: Faces $F_j, F_j'$ of octahedron on the left dual to vertices $\alpha_j, \alpha_j'$ of cube on the right. The face of the cube in gray color corresponds to the dot ``$\bullet$" in the octahedron.
\end{center}
\end{figure}
\noindent 
Consider the rational simplicial decomposition
$\Xi^*$ of
$\triangle$, which is a refinement of $\Xi$ by
adding the center 
$$
c := \frac{1}{4} \sum_{j=1}^4 e^j \ 
$$
as a vertex with the barycentric decomposition
of $\Diamond$ in $\Xi$. Note that $c \not\in N$
and $2c \in N$. See [Fig. 6].
\begin{figure}[ht]
\begin{center}
\mbox{\psfig{figure=tatrhilb.ps}}

\vskip6pt
[Fig. 6]  The rational simplicial decomposition $\Xi^*$ of $\bigtriangleup$ for $n=4, r=1$.  
\end{center}
\end{figure}

\begin{theorem}\label{th:A1(4)}
For $G= A_1(4)$, we have 
$$
{\rm Hilb}^{A_1(4)}(\CZ^4) = X_{A_1(4)} \simeq 
X_{\Xi^*} \ , 
$$
which is non-singular with the canonical bundle 
$\omega = {\cal O}_{X_{\Xi^*}} (E)$, where $E$ is
an irreducible divisor isomorphic to the triple
product of
$\PZ^1$, 
\begin{eqnarray}
E = \PZ^1 \times \PZ^1 \times \PZ^1 \ . 
\label{eq:E}
\end{eqnarray}
Furthermore for $\{ i, j, k \} = \{1,2,3
\}$, the normal bundle of $E$ when
restricted on the $\PZ^1$-fiber, $\PZ^1_k$,
for the projection on $\PZ^1\times \PZ^1$ via the
$(i, j)$-th factor,
\begin{eqnarray}
p_k : E \longrightarrow \PZ^1 \times \PZ^1 \
, 
\label{eq:proj}
\end{eqnarray}
is the
$(-1)$-hyperplane bundle:
\begin{eqnarray}
{\cal O}_{X_{\Xi^*}}(E) \otimes {\cal
O}_{\PZ^1_k}
\simeq  {\cal O}_{\PZ^1}(-1)  \ . 
\label{eq:f(-1)}
\end{eqnarray}
\end{theorem}
{\it Proof.} By Lemma \ref{lem:SingXi} and 
Remark (1) after that, one can see the smoothness
of $X_{\Xi^*}$ on the affine chart corresponding to
$\triangle_j$, also its relation with ${\rm
Hilb}^G(\CZ^4)$. For the rest of simplexes, 
the octahedron $\Diamond$ of $\Xi$ is decomposed
into eight simplexes corresponding to the faces
$F_j, F_j'$ of $\Diamond$. Denote $C_j$ ($C_j'$)
the simplex of $\Xi^*$ spanned by $c$ and $F_j$
($F_j'$ respectively), and $x_{C_j}, x_{C_j'}$
the elements in $X_{\Xi^*}$ of the corresponding 
$T$-orbit. First we show that for $x=x_{C_j},
x_{C_j'}$, ${\cal F}_{X_{\Xi^*, x}}$ is a regular
$G$-module.
It is easy to see that
the vertices of $F_j$ together with $2c$ form a
integral basis of
$N$, the same for the vertices of $F_j'$. For the
convenience of notation, we can set $j=1$ without loss 
of generality. Then we
have the integral basis of $M$ for the cones, dual to  
$C_1, C_1'$ as follows:
$$
\begin{array}{ll}
{\rm cone}(C_1)^* : &(2c, v^{1,2}, v^{1,3},
v^{1,4})^{-1} = 
\left( \begin{array}{cccc}
-1 &   1 & 1 & 1  \\
1 & 1 & -1 & -1 \\
1 & -1 & 1 & -1 \\
1 & -1 & -1 & 1
\end{array} \right) \ , \\
{\rm cone}(C_1')^* : &
(2c, v^{2,3}, v^{2,4}, v^{3,4})^{-1} = 
\left( \begin{array}{cccc}
2 &   0 & 0 & 0  \\
-1 & 1 & 1 & -1 \\
-1 & 1 & -1 & 1 \\
-1 & -1 & 1 & 1
\end{array} \right) \ .
\end{array}
$$
Therefore, the following 4 functions form a
smooth coordinate of $X_{\Xi^*}$ near $x_{C_j}$ for $j=1$,
$$
U_1= \frac{Z_2Z_3Z_4}{Z_1} , \ \ \
U_2= \frac{Z_1Z_2}{Z_3Z_4}, \ \ \
U_3= \frac{Z_1Z_3}{Z_2Z_4}, \ \ \ 
U_4= \frac{Z_1Z_4}{Z_2Z_3} \ ,
$$
and one has
$$
I (x_{C_1}) = < Z_2Z_3Z_4, Z_1Z_2, Z_1Z_3, Z_1Z_4
> + I_G \ .
$$
Similarly the coordinates near $x_{C_j'}$ for $j=1$ are
given by
$$
U_1'= Z_1^2 , \ \ \ 
U_2'= \frac{Z_2Z_3}{Z_1Z_4}, \ \ \ 
U_3'= \frac{Z_2Z_4}{Z_1Z_3}, \ \ \ 
U_4'= \frac{Z_3Z_4}{Z_1Z_2} \ ,
$$
and we have
$$
I (x_{C_1'}) = <Z_2Z_3, Z_2Z_4, Z_3Z_4 > + I_G \ .
$$
It is easy to see that the $G$-modules, $\CZ[Z]/I
(x_{C_1}), 
\CZ[Z]/I (x_{C_1})$, are both equivalent to the
regular representation. Therefore the ideals
$I(x)$ for $x=x_{\triangle_j}, x_{C_j}, x_{C'_j}$
(
$1 \leq j \leq 4$), give rise to distinct elements in
${\rm Hilb}^{A_1(4)}(\CZ^4)$. In fact, one can
show that $X_{\Xi^*} = {\rm Hilb}^{A_1(4)}(\CZ^4)$ 
( for the details, see \cite{CR}). By (\ref{K}), the canonical bundle of
$X_{\Xi^*}$ is given by
$$
\omega_{X_{\Xi^*}} = {\cal O}_{X_{\Xi^*}}(E)
$$
where $E$ is the toric divisor $D_c$. It is
known that $E$ is a 3-dimensional complete toric
variety arisen from the star of $c$ in $\Xi^*$,
which is given by the octahedron in [Fig. 5]; 
in fact, the cube in [Fig. 5] represents the toric orbits' structure.
Therefore $E$ is isomorphic to the triple product
of $\PZ^1$ as in (\ref{eq:E}). The conclusion of the
normal bundle of $E$ restricting on each
$\PZ^1$-fiber will follow from techniques in
toric geometry. For example, for fibers over
the projection of $E$ onto the $(\PZ^1)^2$
corresponding to the 2-convex set spanned
by $v^{1,2}, v^{1,3}, v^{3,4}, v^{2,4}$, one can
perform the computation as follows.
Let $(U_1, U_2, U_3, U_4)$ be
the local coordinates near $x_{C_4'}$ dual to the
$N$-basis 
$(2c, v^{1,2}, v^{1,3},  v^{2,3})$,
similarly the local coordinate $(W_1, W_2, W_3,
W_4)$ near $x_{C_1}$ dual to $(2c, v^{1,2},
v^{1,3}, v^{1,4})$. By $2c= v^{1,4}+v^{2,3}$,
one has the relations, 
$$
U_1=W_1W_4 \ , \ \ U_4= W_4^{-1} \ , \ \ \ U_j=W_j \ (j
=2,3) \ .
$$
This shows that the restriction of the normal
bundle of $E$ on a fiber $\PZ^1$ over $(U_2,
U_3)$-plane is the $(-1)$-hyperplane bundle.
\hfill $\Box$ \par \vspace{.2in} \noindent
The sheaf ${\cal F}_{X_{\Xi^*}}$ for $X_{\Xi^*}$ 
in Theorem
\ref{th:A1(4)} is a vector bundle with the regular
$G$-module on each fiber. The local
frame of the vector bundle is provided by
the structure of $\CZ[Z]/I(x)$ for $x$ being the
zero-dimensional toric orbit of $X_{\Xi^*}$. One
can have a pictorial realization of monomial
basis of these $G$-representations as follows. We
start with the element $x_{\triangle_1}$,
and the identification,
$\CZ[Z]/I(x_{\triangle_1})=
I(x_{\triangle_1})^\bot$. The eigen-basis of the
$G$-module $I(x_{\triangle_1})^\bot$ is given by
monomials in the diagram of [Fig. 7].

\begin{figure}[ht]
\begin{center}
\mbox{\psfig{figure=4splot00.ps}}

[Fig. 7] The monomial basis of the
$G$-module $I(x_{\triangle_1})^\bot (\simeq \CZ[Z]/I(x_{\triangle_1}))$ in the $Z_2$-$Z_3$-$Z_4$ space. 
\end{center}
\end{figure}
\noindent
By the fact that the $\rho$-eigenspace of
$I(o)^\bot$ for the element $o \in S_G$ has
the dimension 2 for a non-trivial character
$\rho$, one can present the data of the regular 
$G$-module $\CZ[Z]/I(x)$ for another element $x$
by indicating the ones to replaced in the [Fig. 7],
which will be marked by $\times$. All of the 
replacement are list in [Fig. 8-10].  

\begin{figure}[ht]
\begin{center}
\mbox{\psfig{figure=4splot09.ps}}
\mbox{\psfig{figure=4splot10.ps}}
\mbox{\psfig{figure=4splot11.ps}}

{\vskip 10pt}
[Fig. 8] The corresponding $I^\bot$-graph for the simplex $\Delta_2$, $\Delta_3$ and $\Delta_4$ in $X_{\Xi^*}$s. An ``$\bullet$" means a monomial in $I(x_{\triangle_i})^\bot$, while an ``$\times$" means one in $I(x_{\triangle_i})$.
\end{center}
\end{figure}

\begin{figure}[ht]
\begin{center}
\mbox{\psfig{figure=4splot01.ps}}
\mbox{\psfig{figure=4splot05.ps}}
\mbox{\psfig{figure=4splot06.ps}}
\mbox{\psfig{figure=4splot07.ps}}

{\vskip 10pt}
[Fig. 9] The corresponding $I^\bot$-graph for the simplex $C_1$, $C_2$, $C_3$ and $C_4$
\end{center}
\end{figure}
\begin{figure}[ht]
\begin{center}
\mbox{\psfig{figure=4splot08.ps}}
\mbox{\psfig{figure=4splot04.ps}}
\mbox{\psfig{figure=4splot03.ps}}
\mbox{\psfig{figure=4splot02.ps}}

{\vskip 10pt}
[Fig. 10] The corresponding $I^\bot$-graph for the simplex $C_1^\prime$, $C_2^\prime$, $C_3^\prime$ and $C_4^\prime$
\end{center}
\end{figure}

By the standard blowing-down
criterion of an exceptional divisor, the property
(\ref{eq:f(-1)}) ensures the existence of a
smooth $4$-fold $(X_{\Xi^*})_k$ by blowing down the
family of
$\PZ^1$s along the projection $p_k$
(\ref{eq:proj}) for each $k$. In fact,
$(X_{\Xi^*})_k$ is also a toric variety
$X_{\Xi_k}$ where $\Xi_k$ is the refinement of 
$\Xi$ by adding the segment connecting
$v^{k,4}$ and $v^{i, j}$ to divide the central
polygon $\Diamond$ into $4$ simplexes, where $\{i,
j, k \}= \{1,2,3 \}$. Each $X_{\Xi_k}$ is a
crepant resolution of $X_\Xi (= S_{A_1(4)})$. 
We have the 
relation of refinements: $
\Xi \prec \Xi_k \prec \Xi^* $ for $k=1, 2, 3$.
The polyhedral decomposition in the
central part  
$\Diamond$ appeared in the refinement relation are
denoted by 
$$
\Diamond \prec \Diamond_k \prec \Diamond^* \ , \ \
\ \ k=1, 2, 3 ,
$$
of which the pictorial realization is given in [Fig. 11]. 
\begin{figure}[ht]
\begin{center}
\mbox{\psfig{figure=hilbs.ps}}

\mbox{\psfig{figure=crep1.ps}} {\hskip 50pt}\mbox{\psfig{figure=crep2.ps}}{\hskip 50pt}\mbox{\psfig{figure=crep3.ps}}

\mbox{\psfig{figure=octa.ps}}

[Fig. 11] Toric representation of 4-dimensional flops over a common singular base and dominated by the same 4-fold.    
\end{center}
\end{figure}
The connection of smooth $4$-folds for 
different
$\Diamond_k$ can be regarded as a ``flop"
of 4-folds suggested by the similar procedure in the theory
of 3-dimensional birational geometry. Each one is
a ``small"\footnote{Here the smallness for a resolution means one with the exceptional locus of codimension $\geq 2$.} resolution of a 4-dimensional 
isolated singularity defined by the equation 
(\ref{4sing}).
Hence we have shown the following result.
\begin{theorem} \label{th:flop}
For $G= A_1(4)$, there are crepant resolutions of
$S_G$ obtained by blowing down the divisor $E$ of ${\rm
Hilb}^G(\CZ^4)$ along $(\ref{eq:proj})$ in Theorem
$\ref{th:A1(4)}$. Any two such resolutions differ
by  a ``flop" procedure of $4$-folds. 
\end{theorem}

\section{${\bf A_r(4)} $-Singularity and Conclusion
Remarks}  
For $G=A_r(n), n \geq 4$,
the structure of ${\rm Hilb}^G(\CZ^n)$ and its
relation with possible crepant resolutions of
$S_G$ has been an on-going program under
investigation. We have discussed the simplest case
$A_1(4)$ in Theorem
\ref{th:A1(4)}. A similar conclusion holds also
for $n=4$, but an arbitrary $r$, whose proof
relies on more complicated techniques. The
details  will be given in \cite{CR}.
\begin{theorem}\label{th:A(4)}
The $G$-Hilbert
scheme ${\rm Hilb}^{A_r(4)}(\CZ^4)$ is the 
non-singular toric variety 
$X_{A_r(4)}$  with the
canonical bundle,
$$
\omega_{X_{A_r(4)}} = {\cal O} ( 
\sum_{k=1}^m E_k ) \ , \ \ m =
\frac{r(r+1)(r+2)}{6} \ ,  
$$ 
where $E_k$s are disjoint smooth
exceptional divisors in
$X_{A_r(4)}$, each of which satisfies the
condition
$(\ref{eq:E}), (\ref{eq:f(-1)})$. Associated to a
projection $(\ref{eq:proj})$ for each $E_k$, there
corresponds a toric  crepant resolution
$\widehat{S}_{A_r(4)}$ of
$S_{A_r(4)}$ with  
\begin{eqnarray*}
\chi ( \widehat{S}_{A_r(4)} ) = |A_r(4) |
= (r+1)^3 . 
\end{eqnarray*}
Furthermore, any two such
$\widehat{S}_{A_r(4)}$s differ by flops of
$4$-fold.
\end{theorem}
One can also describe a monomial
basis of the $A_r(4)$-module $\CZ[Z]/I(x)$ for
$x \in {\rm Hilb}^{A_r(4)}(\CZ^4)$, similar to
the one we have given in the previous section for
the case $r=1$.

Another type of hypersurface orbifolds are those
from the quotient singularities of simple groups.
For $n=3$, the well-known examples are
$S_G$ for $G$= the icosahedral group
$I_{60}$, Klein group $H_{168}$, in which cases
a crepant resolution of $S_G$ was explicitly
constructed in \cite{R94}, \cite{M}
respectively. The structure of ${\rm Hilb}^G(\CZ^3)$ has recently been discussed in \cite{GNS}. Even though the crepant and smooth property of the group orbit Hilbert
schemes for dimension 3 is known by \cite{BKR}, a clear quantitative and qualitative relation would still be interesting for the possible study of some other simple groups $G$ in higher dimension. A such program is under consideration with initial progress being made.

Even for the abelian group $G$ in the dimension
$n=3$, the conclusion on the trivial canonical
bundle of ${\rm Hilb}^G(\CZ^3)$ would raise a
subtle question in  the mirror problem of
Calabi-Yau 3-folds in string theory. As an example,
a standard well known one is the Fermat
quintic in
$\PZ^4$ with the special marginal deformation
family:
$$
X : \ \ \sum_{j=1}^5 z_i^5 + \lambda
z_1z_2z_3z_4z_5 
= 0 \ , \ \ \lambda \in \CZ \ .
$$
With the maximal diagonal group $SD$ of
$z_i$s preserving the family $X$, the mirror
$X^*$ is constructed by ``the" crepant resolution
of $X/SD$, $X^*= \widehat{X/SD}$ (see e.g., 
\cite{GP, R91}), by which the roles of $H^{1,1},
H^{2,1}$ are interchangeable in the ``quantum"
sense. When working on the
one-dimensional space $H^{1,1}(X)
\sim H^{2,1}(X^*)$, the choice of
crepant resolution  $\widehat{X/SD}$ makes no
difference on the conclusion. While on the part of
$H^{2,1}(X) \sim H^{1,1}(X^*)$,
it has been known that many topological
invariants, like Euler characteristic, Hodge
numbers, elliptic genus, are independent of the
choices of crepant resolutions, hence one obtains the same
invariants for different choices of crepant resolutions as the model for $X^*$. However, the topological triple
intersection of cohomologies does differ for two
crepant resolutions (see, e.g.,
\cite{R93}), hence the  choice of
crepant resolution as the mirror 
$X^*=\widehat{X/SD}$ will lead to the different
effect on the topological cubic form of
$H^{1,1}(X^*)$, upon which as the ``classical" level,
the quantum triple  product of the physical mirror 
theory would be built (see, e.g., articles in
\cite{Y}). The question of the
"good" model for $X^*$ has rarely
been raised in the past, partly due to the lack of
mathematical knowledge on the issue. However, 
with the
$G$-Hilbert scheme now given in
Sect. 3, 4 as the mirror $X^*$, it
seems to have left some
fundamental  open problems on its formalism of
mirror Calabi-Yau spaces and the
question of the arbitrariness of the choice
of crepant resolutions remains a mathematical
question to be completely understood concerning its applicable 
physical theory. 

For the role of $G$-Hilbert scheme in the study
of crepant resolution of $S_G$, the conclusion we
have obtained for
$G= A_r(4)$ has indicated that  ${\rm
Hilb}^G(\CZ^n)$ couldn't be a crepant resolution
of $S_G$ in general when the dimension $n$ is
greater than 3. Nevertheless the structure of ${\rm
Hilb}^G(\CZ^n)$ is worthwhile for further study on its
own right due to the interplay of geometry and group
representations. Its understanding could still lead to
the construction of crepant resolutions
of $S_G$ in case such one does exist. It would be
a promising direction of the geometrical study
of orbifolds.

\section{Acknowledgements}
We wish to thank I. Nakamura for making us aware of his
work \cite{GNS, Na}. This work is supported in
part by the
NSC of Taiwan under grant No.89-2115-M-001-012.


\begin{thebibliography}{99}
\bibitem{BKR} T. Bridgeland, A. King and M. Reid, The McKay correspondence as an equivalence of derived categories, J. Amer. Math. Soc. 14(2001), no. 3, 535-554.

\bibitem{CR} L. Chiang and S. S. Roan, On hypersurface quotient singularity of dimension 3 and 4, (math.AG/0011151). 

\bibitem{CLO} D. Cox, J. Little and D. O'Shea,
Ideals, Varieties, and Algorithms,
Springer-Verlag, New York-Berlin-Heidelburg, 1992.

\bibitem{Dan} V. I. Danilov, The geometry of toric varieties, Russ. Math. 
Surveys 33:2 (1978) 97-154.

\bibitem{DHVW} L. Dixon, J. Harvey, C. Vafa and E. Witten, Strings on
orbifolds, Nucl. Phys. B 261 (1985) 678-686, Strings on orbifolds II,
Nucl. Phys. B 274 (1986) 285-314.

\bibitem{GNS} Y. Gomi, I. Nakamura and K. Shinoda, Hilbert schemes of $G$-orbits in dimension three, Asian J. Math. 4(1) (2000) 51-70. 

\bibitem{GV} G. Gonazles-Springer and J. L.
Verdier, Construction g\'{e}ometrique de la
correspondence de McKay, Ann. Sci. \'{E}cole
Norm. Sup. 16 (1983) 409-449.

\bibitem{GP} B. R. Greene and M. R. Plesser,
Duality in Calabi-Yau moduli space, Nucl. Phys. B
338 (1990) 15-37.

\bibitem{HH} F. Hirzebruch and T. H\"{o}fer,
On the Euler number  of an orbifold, Math. Ann.
286 (1990) 255-260.

\bibitem{INj} Y. Ito and H. Nakajima, McKay
correspondence and Hilbert schemes in dimension
three, Topology 39 (2000) 1155-1191, (math.AG/9803120).

\bibitem{IN} Y. Ito and I. Nakamura, McKay
correspondence and Hilbert schemes, Proc.
Japan Acad. 72 (1996) 135-138, Hilbert  schemes
and simple  singularities, New trends in
Algebraic Geometry (Proc. of the July 1996
Warwick European Alg. Geom. Conf.) Cambridge
Univ. Press (1999) 151-233.

I. Nakamura, Hilbert scheme and simple singularities $E_6, E_7$ and $E_8$, Hokkaido Univ. Preprint Series in Math. No. 362 (1996).

\bibitem{KKMS} G. Kempt, F. Knudson, D. Mumford
and B. Saint-Donat, Toroidal embedding 1, Lecture
Notes in Math., Vol. 339, Springer-Verlag, New York,
1973.

\bibitem{Kl} F. Klein, Gesammelte Mathematische
Abhandlungen Bd. II,  Springer-Verlag 1922
(reprint 1973).  

\bibitem{M} D. Markushevich, Resolution of 
${\bf C}^3/H_{168}$, Math. Ann. 308(1997) 279-289.

\bibitem{Mc} J. McKay, Graphs, singularities and finite groups, Proc. Symp. Pure Math. 37(1980) 183-186.

\bibitem{Mo} S. Mori, Birational classification of algebraic
threefolds,  Internat. Cong. Math. Kyoto (1990) 235-248.


\bibitem{Na} I. Nakamura, Hilbert schemes of
abelian group orbits, J. Alg. Geom. 10(2001)757-779.

\bibitem{O} T. Oda, Lectures on torus embeddings and applications, Tata 
Inst. of Fund. Research 58, Springer-Verlag (1978).


\bibitem{R89} S. S. Roan, On the generalization 
of Kummer surfaces, J. Diff. Geom. 30(1989) 523-537.

\bibitem{R90} S. S. Roan, On the Calabi-Yau
orbifolds in weighted projective spaces, 
Internat. J. Math.  1(2)(1990) 211-232.

\bibitem{R91} S. S. Roan, The mirror of 
Calabi-Yau orbifold, Internat. J. Math. 2(4)(1991) 439-455.

\bibitem{R93} S. S. Roan, Topological coupling
of  Calabi-Yau orbifold,  J. Group Theory in
Phys. 1(1993) 83- 103.

\bibitem{R94}  S. S. Roan, On $c_1 = 0$ resolution
of quotient singularity,  Internat. J. Math. 5(4)(1994) 523-536.

\bibitem{Rtop} S. S. Roan,  Minimal resolutions
of Gorenstein orbifolds in dimension three,
Topology 35(1996) 489-508.

\bibitem{Y} S. T. Yau  ( ed. ), 
Essays on mirror manifolds,  Internat. Press, 1992.

\end{thebibliography}
\end{document}